\magnification = 1200
\showboxdepth=0 \showboxbreadth=0
 %no black boxes

\baselineskip 15pt
\parskip3pt
%\voffset = -0.5in
%\hoffset = -.05truein
\def\qed{\hfill\vrule height6pt width6pt depth0pt}

\def\ss{\smallskip}

\def\bs{\bigskip}

\def\cl{\centerline}

\def\nind{\noindent}

\def\ref#1#2{\nind\hangindent.5in\hbox to .5in{#1\hfill}#2}
\def\reff#1#2{\nind\hangindent.8in\hbox to .8in{\bf #1\hfill}#2\par}
\def\refd#1#2{\nind\hangindent.8in\hbox to .8in{\bf #1\hfill {\rm
--}}#2\par}
\def\pmb#1{\setbox0=\hbox{#1}
\kern-0.025em\copy0\kern-\wd0
\kern.05em\copy0\kern-\wd0
\kern-.025em\raise.0433em\box0}

\def\frac#1#2{{#1\over#2}}

\def\text#1{\rm#1}

\outer\def\stmnt  #1. #2\par{\medbreak
\noindent{\bf#1.\enspace}{\sl#2}\par
\ifdim\lastskip<\medskipamount \removelastskip\penalty55\medskip\fi}

\def\newline{\hfill\break}
\def\:{\,:\,}

\def\({\left(}                   
\def\){\right)}

\def\[{\left[}                   
\def\]{\right]}

\def\lan{\langle}
\def\ran{\rangle}

\def\ci{\subset}

\def\fy{\infty}

\def\lam{\lambda}

\def\om{\omega}

\def\Om{\Omega}

\cl{\bf Viscosity convex functions on Carnot groups}
\bs
\cl{Changyou Wang}
\ss
\cl{Department of Mathematics, \ University of Kentucky}
\cl{Lexington, KY 40506}  
\bs
\cl{\bf Abstract} 
\ss
{\it We prove that any locally bounded from below, upper semicontinuous  
v-convex function in any Carnot group is  h-convex. }
\bs
\nind{\S1}. Introduction

Convex functions have played very important roles in PDEs,
especially fully nonlinear elliptic PDEs in Euclidean spaces
(see Caffarelli-Cabr\'e [CC] and Crandall-Ishii-Lions [CIL]).
Motivated by this fact and the aim to develop an
intrinsic theory of subelliptic fully nonlinear PDEs on
Carnot groups, there have been works towards
the theory of convex functions on Heisenberg groups 
by Lu-Manfredi-Stroffolini [LMS], and on general Carnot
groups by Danielli-Garofalo-Nhieu [DGN].

Lu-Manfredi-Stroffolini [LMS] have extended the concept
of {\it convex in the viscosity sense} (or v-convex)
from the Euclidean space to the sub-Riemannian setting of Heisenberg
groups. Using the uniqueness theorem
on viscosity solutions of the subelliptic $\fy$-laplacian
by Bieske [B], they showed that any upper semicontinuous 
v-convex function on any Heisenberg group is locally Lipschitz continuous.

A geometric approach of convexity on Carnot groups is given
by Danielli-Carofalo-Nhieu [DGN], where they have introduced
the notation of {\it horizontally convex} (or h-convex)
functions. One of the theorems of [DGN] implies that any
locally bounded h-convex function is locally Lipschitz
continuous. 

It is known by [LMS] that any upper semicontinuous h-convex
function $u$ on $\bf G$ is v-convex, and the converse
is also true if, in addition, $u\in \Gamma^2(\bf G)$
(the horizontal $C^2$ space). Here we are interested in
whether the converse remains true under minimal regularity
assumptions.

In order to state our theorem, we first recall the basic properties of Carnot groups.
A simply connected Lie group $\bf G$ is called a Carnot group
of step $r\ge 1$, if its Lie algebra ${\it g}$ admits a vector space
decomposition in $r$ layers ${\it g}=V_1+V_2+\cdots +V_r$ such that
(i) ${\it g}$ is stratified, i.e., $[V_1, V_j]=V_{j+1}, j=1,\cdots, r-1$,
and (ii) ${\it g}$ is $r$-nilpotent, i.e. $[V_j, V_r]=0, j=1, \cdots, r$.  
We call $V_1$ the {\it horizontal} layer and $V_j, j=2,\cdots, r$
the {\it vertical} layers.
We choose an inner product $\lan\cdot,\cdot\ran$ on {\it g}
such that $V_j's$ are mutually orthogonal for $1\le j\le r$.
Let $\{X_{j, 1}, \cdots, X_{j,m_j}\}$ denote
a fixed orthonormal basis of $V_j$ for $1\le j\le r$,
where $m_j=\hbox{dim}(V_j)$ is the dimension of $V_j$.
From now on, we also denote $m=\hbox{dim}(V_1)$ as
the dimension of the horizontal layer and set
$X_i=X_{1,i}$ for $1\le i\le m$. It is well-known (see [FS])
that the exponential map ${\it exp}: {\it g}\equiv R^n\to G$ is a global
diffeomorphism so that there exists an exponential coordinate
system on $\bf G$ with $n=m+\sum_{i=2}^r m_i$ as its topological dimension. 
More precisely, any $p\in \bf G$ has a coordinate 
$((p_1, \cdots, p_m),(p_{2,1},\cdots, p_{2,m_2}), \cdots,
(p_{r,1},\cdots, p_{r,m_r}))$ such that
$$p={\it exp}(\xi_1(p)+\cdots \xi_r(p)), \hbox{ with }
\xi_1(p)=\sum_{l=1}^m p_lX_l, \xi_i(p)=\sum_{j=1}^{m_i}p_{i,j}X_{i,j}, 2\le i\le r.
\eqno(1.1)$$
The exponential map can induce a homogeneous pseudo-norm $N_{\bf G}$ on $\bf G$ 
in the following way (see [FS]). 
$$N_{\bf G}(p):=(\sum_{i=1}^r|\xi_i(p)|^{2r!\over i})^{1\over 2r!}, \hbox{ if }
p={\it exp}(\xi_1(p)+\cdots \xi_r(p)),\eqno(1.2)$$
where $|\xi_1(p)|=(\sum_{l=1}^m p_l^2)^{1\over 2}$ and
$|\xi_i(p)|=(\sum_{j=1}^{m_i} p_{i,j}^2)^{1\over 2}(2\le i\le r)$.
Moreover, $N_{\bf G}$ yields a pseudo-distance on $\bf G$ as follows.
$$d_{\bf G}(p, q):=N_{\bf G}(p^{-1}\cdot q), \ \forall p,\ q\in\bf G,\eqno(1.3)$$
where ${\cdot}$ is the group multiplication of ${\bf G}$ and
$p^{-1}$ is the inverse of $p$. It is easy to see
that $d_{\bf G}$ satisfies the invariance property
$$d_{\bf G}(z\cdot x,z\cdot y)=d_{\bf G}(x,y),  \ \ \forall x, \ y, \ z\in G,\eqno(1.4)$$
and is of homogeneous of degree one, i.e.
$$d_{\bf G}(\delta_\lam(p), \delta_\lam(q))=\lam d_{\bf G}(p,q),\ \forall \lam>0,
\ \forall p, \ q\in\bf G \eqno(1.5)$$
where $\delta_\lam(p)=\lam\xi_1(p)+\sum_{i=2}^r \lam ^i\xi_i(p)$ is the nonisotropic
dilations on $\bf G$. 

We need some notations. For $u:\bf G\to R$,
let $\nabla_h u:=(X_1 u, \cdots, X_m u)$ denote the horizontal gradient of $u$,
$\nabla^2_h u:=({X_i X_j +X_jX_i\over 2}u)_{1\le i, j\le m}$ denote the
horizontal hessian of $u$,
and use $\nabla u$ and $\nabla^2 u$ denote the full gradient and hessian of $u$
respectively.  For any $p\in{\bf G}$, 
let $H_p({\bf G})=\hbox{span}\{X_1(p),\cdots, X_m(p)\}$ denote the horizontal
tangent plane of ${\bf G}$ at $p$. 

We now recall the definition of {\it horizontal convexity} introduced by [DGN]
(\S5, Definition 5.5), see also [LMS](\S4, Definition 4.1).
\ss
\nind{\bf Definition 1.1}. For a domain $\Om\ci{\bf G}$, a function $u:\Om\to R$
is called to be {\it horizontally convex} (h-convex) if
for any $p\in\Om$ and $q\in H_p({\bf G})\cap\Om$, $u|_{[p,q]}$ is convex,
where $[p,q]$ denotes the line segment joining $p$ and $q$.

The v-convexity has been introduced by [LMS](\S3, Definition 3.1), see
[CIL] for the general theory of viscosity solutions. More precisely,
\ss
\nind{\bf Definition 1.2}. For a domain $\Om\ci{\bf G}$, a upper semicontinuous
function $u:\Om\to R$ is  convex in the viscosity sense (v-convex) if for any
vector $\xi\in R^m$, $\xi^T\nabla^2_h u\xi\ge 0$ in the sense of viscosity, i.e.
for any $p\in\Om$ and any $\phi\in C^2(\Om)$ touching $u$ from above at $p$, there holds
$$\xi^T\nabla^2_h \phi(p)\xi=\sum_{i, j=1}^m \xi_i\xi_j X_iX_j\phi(p)\ge 0.
\eqno(1.6)$$

We are ready to state
\ss
\nind{\bf Theorem A}. {\it Any locally bounded from below, upper semicontinuous v-convex
functions on a Carnot group $\bf G$ are h-convex.}

\ss
\nind{\bf Remark 1.3}. Balogh-Rickly [BR] has recently proved
theorem A for Heisenberg groups without the assumption
of local lower boundedness. While preparing this paper,
J. Manfredi has also informed me that Juutinen-Manfredi [JM] are
able to prove theorem A by a different method.

Our idea to prove theorem A is based on the sup-convolution construction
(see \S2 below) on Carnot groups, which was developed in an earlier 
paper by Wang [W] and was employed to prove the uniqueness for continuous
viscosity solutions to the subelliptic $\fy$-laplacian equations on
any Carnot group $\bf G$. Roughly speaking, the sup-convolution
of a v-convex function is not only v-convex but also semiconvex in
the Euclidean sense. The reader can consult with Jensen-Lions-Sougannidis [JLS]
for the sup-convolution in the Euclidean space.

As a byproduct of the proof of theorem A, we also obtain the following characterization
of continuous v-convex functions, analogous to that of convex functions
on the Euclidean space in the viscosity sense (cf. [LMS]\S2).
\ss
\nind{\bf Corollary B}.  {\it For any bounded domain $\Om\ci\bf G$
and $u\in C({\Om})$. Then the following
statements are equivalent:

(a) $u$ is v-convex on $\Om$.

(b) For any subdomain $\tilde{\Om}\ci\ci\Om$,
there exist $\{u_k\}\ci \Gamma^2({\tilde\Om})$ such that $u_k$'s are v-convex on ${\tilde\Om}$
and $u_k\rightarrow u$ uniformly on ${\tilde\Om}$.

(c) For any subdomain $\tilde{\Om}\ci\ci\Om$,
there exist $\{u_k\}\ci \Gamma^2({\tilde\Om})$ such that $u_k$'s are h-convex on ${\tilde\Om}$
and $u_k\rightarrow u$ uniformly on ${\tilde\Om}$.}

\ss
The note is written as follows. In \S2, we outline the sup-convolution construction.
In \S3, we prove theorem A. 
 
\bs
\nind\S2 The construction of sup-convolutios on $\bf G$
\ss
In this section, we outline the construction of sup-convolutions
on $\bf G$, which was developed by Wang [W] on Carnot groups
and by [JLS] on Euclidean spaces, to show that the sup-convolution of
a v-convex function is v-convex.

Let $\Om\ci\bf G$ be a bounded domain 
and $d_{\bf G}(\cdot,\cdot)$ be the smooth gauge pseudo-distance
defined by (1.3). For any $\epsilon>0$, define
$$\Om_\epsilon=\{x\in \Om: \inf_{y\in {\bf G}\setminus\Om} d_{\bf G}(x^{-1},y^{-1})
\ge\epsilon\}.$$
\ss
\nind{\bf Definition 2.1}. For $\epsilon>0$ and $u:\bar\Om\to R$ a upper semicontinous
and bounded from below function, the
sup-convolution $u_\epsilon$ of $u$ is defined by
$$u^\epsilon(x)=\sup_{y\in\bar\Om}(u(y)-{1\over 2\epsilon}d_{\bf G}(x^{-1},y^{-1})^{2r!}), 
\ \forall x\in\Om. \eqno(2.1)$$

For $p\in \bf G$, we also denote $\|p\|_{E}:=(\sum_{i=1}^m |\xi_i(p)|^2)^{1\over 2}$ 
as the euclidean norm of $p$. We recall
\ss
\nind{\bf Definition 2.2}. A upper semicontinuous function $u:\bar\Om\to R$ 
is called semiconvex, if
there is a constant $C>0$ such that $u(p)+C\|p\|_E^2:\bar\Om\to R$ is convex
in the Euclidean sense. This roughly means that if $u\in C^2(\bar\Om)$ and
$\nabla ^2u(p)+C$ is positive semidefinite for $p\in\Om$, then $u$ is semiconvex.

Now we have the following generalized version of [JLS], 
which can be found in [W](\S3, Proposition 3.3).
\ss
\nind{\bf Proposition 2.3}. {\it Suppose that $u:\bar\Om\to R$ is
upper semicontinuous and  bounded from below, and let
$R_0=\|u\|_{L^\fy(\Om)}>0$. Then
$u_\epsilon>0$ satisfies:

(1) $u^\epsilon$ is locally Lipschitz continuous in $\Om$ with respect to $d_{\bf G}$,
and is semiconvex.

(2) $\{u^\epsilon\}$ is monotonically nondecreasing w.r.t. $\epsilon$ 
and converges to $u$ in $\Om$.

(3) if $u$ is a v-convex function in $\Om$,
then $u^\epsilon$ is a v-convex function in $\Om_{(2R_0+1)\epsilon}$.

(4) if, in addition, $u\in C(\bar\Om)$, then $u_\epsilon\rightarrow u$ uniformly on $\Om$.}
\ss
\nind{\bf Proof}. For $\Om\ci\bf G$ is bounded, 
the formula (1.3) for $d_{\bf G}$ implies
$$C(\Om,d)\equiv \|\nabla^2_x(d_{\bf G}(x^{-1},y^{-1})^{2r!})
\|_{L^\fy(\bar\Om\times\bar\Om)}<\fy.$$
Therefore, for any $y\in\bar\Om$, the full hessian of
$$\tilde {u^\epsilon}(x,y):=u(y)-{1\over 2\epsilon}
d_{\bf G}(x^{-1},y^{-1})^{2r!}+{C(\Om,d)\over 2\epsilon}\|x\|_{E}^2,
\ \forall x\in\Om,$$
is positive semidefinite so that ${\tilde{u^\epsilon}}$ is convex.
Note that the superum for a family of convex functions is still convex.
Hence we have
$$\sup_{y\in\bar\Om}{\tilde {u^\epsilon}}(x,y)=u^\epsilon(x)+{C(\Om,d)\over 2\epsilon}\|x\|_{E}^2, 
\ \ x\in\Om$$
is convex. Hence $u^\epsilon$ is semiconvex in $\Om$. It is well-known that semiconvex functions are
locally Lipschitz continuous with respect to the euclidean metric. 
Therefore $u^\epsilon$ is locally Lipschitz continuous with respect to $d_{\bf G}$.
This gives (1). 

For any $\epsilon_1<\epsilon_2$, 
it is easy to see that $u^{\epsilon_1}(x)\le u^{\epsilon_2}(x)$ for any $x\in\Om$
so that $\{u^\epsilon\}$ is monotonically nondecreasing with respect to $\epsilon$.
For any $x\in\Om$, there exists a $x_\epsilon\in\bar\Om$ such that
$$u(x_\epsilon)-{1\over 2\epsilon}d_{\bf G}(x^{-1}, x_\epsilon^{-1})^{2r!}=u^\epsilon(x). \eqno(2.2)$$
This implies
$${1\over 2\epsilon}d_{\bf G}(x^{-1}, x_\epsilon^{-1})^{2r!}\le u(x_\epsilon)-u^\epsilon(x). \eqno(2.3)$$
Since
$$u(x)\le u^\epsilon(x)\le \sup_{y\in \bar\Om}u(y), \ \forall x\in\Om, \eqno(2.4)$$
we have
$${1\over 2\epsilon}d_{\bf G}(x^{-1}, x_\epsilon^{-1})^{2r!}\le 2R_0, \forall x\in\Om. \eqno(2.5)$$
Hence $x_\epsilon\rightarrow x$. Moreover, since
$$u_{\epsilon\over 2}(x)\ge u(x_\epsilon)-{1\over\epsilon}d_{\bf G}(x^{-1}, x_\epsilon^{-1})^{2r!}
=u_\epsilon(x)-{1\over 2\epsilon}d_{\bf G}(x^{-1}, x_\epsilon^{-1})^{2r!},$$
we have
$$\lim_{\epsilon\rightarrow 0}{1\over\epsilon}d_{\bf G}(x^{-1}, x_\epsilon^{-1})^{2r!}=0$$
Therefore, the upper semicontinuity of $u$ implies that $\lim_{\epsilon\downarrow 0}u^\epsilon(x)=u(x)$
for any $x\in\Om$. This gives (2).

For (3), we first observe that for any $x^0\in \Om_{(1+2R_0)\epsilon}$ 
there exists a $x^0_\epsilon\in \Om$ such that 
$$u_\epsilon(x^0)=u(x^0_\epsilon)-{1\over 2\epsilon}d_{\bf G}((x^0)^{-1}, (x^0_\epsilon)^{-1})^{2r!}.$$
Let $\phi\in C^2(\Om)$ be such that
$$u_\epsilon(x^0)-\phi(x^0)\ge u_\epsilon(x)-\phi(x), \ \ \forall x\in\Om.$$
Then we have
$$u(x^0_\epsilon)-{1\over 2\epsilon}d_{\bf G}((x^0)^{-1},(x^0_\epsilon)^{-1})^{2r!}-\phi(x^0)
\ge u(y)-{1\over 2\epsilon}d_{\bf G}(x^{-1}, y^{-1})^{2r!}-\phi(x), \forall x, y\in\Om.$$
For any $y$ close to $x^0_\epsilon$, note that $x=x^0\cdot (x^0_\epsilon)^{-1}
\cdot y\in \Om$. By substituting it into the above inequality, we obtain
$$u(x^0_\epsilon)-\phi(x^0)\ge u(y)-\phi(x^0\cdot (x^0_\epsilon)^{-1}\cdot y).$$
Set $\bar{\phi}(y)=\phi(x_0\cdot (x^0_\epsilon)^{-1}\cdot y)$, for $y\in\Om$
near $y_0$. Then $\bar{\phi}$ touches $u$ at $y=x^0_\epsilon$ from above
so that  v-convexity of $u$ implies
$$\nabla^2_h\bar{\phi}(x^0_\epsilon)\ge 0.$$
On the other hand, the left-invariance of $\{X_i\}_{i=1}^m$ implies
$$\nabla^2_h\bar{\phi}(y)=(\nabla^2_h\phi)(x^0\cdot (x^0_\epsilon)^{-1}\cdot y).$$
Therefore we have
$$\nabla^2_h\phi(x^0))\ge 0,$$
this implies that $u^\epsilon$ is v-convex in $\Om_{(2R_0+1)\epsilon}$.

For (4), since $u\in C(\bar\Om)$, we can see easily from (2.3)-(2.4) that
$$\eqalignno{|u^\epsilon(x)-u(x)|&\le
|u(x_\epsilon)-u(x)|+{1\over 2\epsilon}d_{\bf G}(x^{-1}, x_\epsilon^{-1})^{2r!}\cr
&\le 2|u(x_\epsilon)-u(x)|\le 2\om(d_{\bf G}(x, x_\epsilon)), &(2.6)\cr}$$
where $\om$ is the modular of continuity of $u$.   On the other hand, (2.5) implies
$$d_{\bf G}(x_\epsilon, x)\le (4\epsilon\|u\|_{C(\bar\Om)})^{1\over 2r!}. \eqno(2.7)$$
Therefore, we have
$$\max_{x\in\Om}|u^\epsilon(x)-u(x)|\le 2\om((4\epsilon\|u\|_{C(\bar\Om)})^{1\over 2r!})\rightarrow 0,
\ \hbox{ as }\epsilon\rightarrow 0. \eqno(2.8)$$
The proof is complete.           \qed

\bs
\nind {\S3}. Proof of theorem A and Corollary B
\ss 
This section is devoted to the proof of theorem A. Note that it suffices to prove theorem A
when restricted to any bounded domain of $\bf G$. 
\ss
\nind{\bf Proof}. For any bounded domian $\Om\ci\bf G$ and $\epsilon>0$, 
let $u^\epsilon:\Om\to R$ be the sup-convolution of $u$ obtained by Proposition 2.3. 
Then (1) and (3) of Proposition 2.3 imply 
that $u^\epsilon$ is both semiconvex in the Euclidean sense and v-convex
in $\Om_{(2R_0+1)\epsilon}$. On the other hand, it follows from the well-known theorem
on convex functions in the Euclidean space (cf. Evans-Gariepy [EG]) that $u^\epsilon$
is twice differential in the Euclidean sense for a.e. $x\in \Om_{(2R_0+1)\epsilon}$ 
and $\nabla^2u^\epsilon\in L^1(\Om_{(2R_0+1)\epsilon})$. In particular, 
the horizontal hessian $\nabla^2_hu^\epsilon (x)$ exists for a.e. $x\in\Om_{(2R_0+1)\epsilon}$ and
$\nabla^2_h u^\epsilon\in L^1(\Om_{(2R_0+1)\epsilon})$. 
Since $u^\epsilon$ is v-convex in $\Om_{(2R_0+1)\epsilon}$,
the standard theory on viscosity solutions (see [CIL]) implies that $\nabla^2_hu^\epsilon(x)\ge 0$ 
is positive semidefinite for a.e. $x\in\Om_{(2R_0+1)\epsilon}$, i.e.
$$\sum_{i, j=1}^m \eta_i\eta_j X_iX_j u^\epsilon(x)\ge 0, \ \forall \eta\in R^m. \eqno(3.1)$$
Now let $\phi\in C_0^\fy(\bf G)$ be nonnegative such that supp$(\phi)\ci B_1(0)$ and
$\int_{\bf G}\phi(p)\,dp=1$. For any small $\delta>0$, consider the mollification
$u^\epsilon_\delta=\phi_\delta* u^\epsilon$ of $u^\epsilon$ defined by
$$\phi_\delta*u^\epsilon(q)=\int_{\bf G}u^\epsilon(p^{-1}\cdot q)\phi_\epsilon(p)\,dp,
\ \forall q\in\Om_{(2R_0+1)\epsilon+\delta},$$ 
where $\phi_t(p)=t^{-n}\phi(\delta_t(p))$ for $t>0$.
Since 
$$\nabla^2_h u^\epsilon_\delta(p)=\int_{\bf G}(\nabla^2_h u^\epsilon)(q^{-1}\cdot p)
\phi_\delta(q)\,dq, \ \ \forall p\in \Om_{(2R_0+1)\epsilon+\delta},$$
and $\nabla^2_h u^\epsilon\in L^1(\Om_{(2R_0+1)\epsilon})$ is positive semidefinite, 
we have 
$u^\epsilon_\delta\in C^\fy(\Om_{(2R_0+1)\epsilon+\delta})$ and $\nabla^2_h u^\epsilon_\delta$
is positive semidefinite everywhere in $\Om_{(2R_0+1)\epsilon+\delta}$.
Therefore [DGN](\S5, Theorem 5.11) or [LMS](\S4, Proposition 4.1) implies 
that $u^\epsilon_\delta$ is h-convex on $\Om_{(2R_0+1)\epsilon+\delta}$, i.e.
$$u^\epsilon_\delta(p)
\le {u^\epsilon_\delta(p\cdot h)+u^\epsilon_\delta(p\cdot h^{-1})\over 2}, 
\ \forall p\in \Om_{(2R_0+1)\epsilon+\delta},\  h\in H_0({\bf G})\cap \Om. \eqno(3.2)$$
Since $u^\epsilon_\delta\rightarrow u^\epsilon$ uniformly on $\Om_{(2R_0+1)\epsilon}$ as
$\delta\downarrow 0$, it follows that, by taking $\delta$ into zero in (3.2), 
$u^\epsilon$ is h-convex on $\Om_{(2R_0+1)\epsilon}$. Since (2) of Proposition 2.3
implies $u^\epsilon\rightarrow u$ on $\Om$ as
$\epsilon\downarrow 0$, we can conclude that $u$ is h-convex on $\Om$. \qed
\ss
\nind{\bf Proof of Corollary B}. It is clear that (b) and (c) are equivalent (see, e.g., [DGN][LMS]).
It is also easy to see that (b) implies (a). To see (a) implies (b), let 
${\tilde\Om}\ci\ci\Om$ be fixed and $u^\epsilon$
be the sup-convolution of $u$ on $\Om$ given by Proposition 2.3, 
and $u^\epsilon_{\delta}$ be the $\delta$-mollifer of $u^\epsilon$ on $\Om$
constructed in the proof of theorem A,  where $\epsilon>0, \delta>0$ are sufficiently small.
Then (4) of Proposition 2.3 implies $u^\epsilon\rightarrow u$ uniformly on $\Om$. Moreover,
the proof of theorem A implies that $u_{\delta}^\epsilon\in \Gamma^2({\tilde\Om})$ is v-convex
and $u_{\delta}^\epsilon\rightarrow u^\epsilon$ uniformly on ${\tilde\Om}$.
Therefore, by the Cauchy diagonal process, we may assume that $u_{\delta}^\epsilon\rightarrow
u$ uniformly on ${\tilde\Om}$ for $\epsilon\rightarrow 0$ and $\delta=\delta(\epsilon)\rightarrow 0$.
This finishes the proof.    \qed

\bs
\bs
\cl{\bf REFERENCES}

\nind{[B]} T. Bieske, {\it On $\infty$-harmonic functions on the Heisenberg group}. 
Comm. Partial Differential Equations 27 (2002), no. 3-4, 727--761.

\nind{[BR]} Z. Balogh, M. Rickly, {\it Regularity of convex functions on Heisenberg groups}.
Preprint, 2003.

\nind{[CC]} X. Cabre, L. Caffarelli, {Fully nonlinear elliptic equations}. AMS colloquium publications
43, AMS, Providence, RI, 1995.

\nind{[CIL]} M. Crandall, H. Ishii, P. L. Lions, {\it User's guide to viscosity solutions
of second order partial differential equations}.
Bull. Amer. Math. Soc. (N.S.) 27 (1992), no. 1, 1--67.

\nind{[DGN]} D. Danielli, N. Garofalo, D. Nhieu, {\it Notations of convexity in Carnot groups}.
Comm. Anal. Geom., to appear.

\nind{[EG]} L. Evans, R. Gariepy, {Measure theory and fine properties of functions}. CRC Press, 1992.

\nind{[FS]} G. Folland, E. Stein, Hardy spaces on homogeneous groups. Mathematical Notes, 28. 
Princeton University Press, Princeton, N.J., 1982. 

\nind{[JLS]} R. Jensen,  P. L. Lions, P. Souganidis, {\it A uniqueness result for viscosity
solutions of second order fully nonlinear partial differential equations}. 
Proc. Amer. Math. Soc. 102 (1988), no. 4, 975--978. 

\nind{[JM]} P. Juutinen,  J. Manfredi. In preparation.

\nind{[LMS]} G. Lu, J. Manfredi, B. Stroffolini, {\it Convex functions on Heisenberg group}.
Calc. Var., to appear.

\nind{[W]} C. Y. Wang, {\it The Aronsson equation for absolute minimizers of $L^\fy$-functionals
associated with vector fields satisfying H\"ormander's condition}. Preprint (2003), avaiable
at http//arXiv: math.AP/0307198.

\end